\documentstyle[12pt,psfig]{amsart}

\newtheorem{th}{Theorem}[section]	
\newtheorem{lem}[th]{Lemma}		
		\newtheorem{prp}[th]{Proposition} 

\theoremstyle{definition}       

		\newtheorem{question}[th]{Question}

\theoremstyle{remark}		

\newcommand{\rar}{\rightarrow} 	
	
   	\newcommand{\bfc}{{\Bbb{C}}}
	
\newcommand{\bfz}{{\Bbb{Z}}}	
\newcommand{\bfr}{{\Bbb{R}}}
	\newcommand{\I}{{\cal{I}}}

\newcommand{\Aut}{{\operatorname{Aut}}}

\newcommand{\red}{{\mbox{\small red}}}

\newcommand{\ord}{{\operatorname{ord}}}

\begin{document}
\large

\title{Equivariant Resolution of Singurlarities in Characteristic 0}
\maketitle

\noindent
\today\ (preliminary Version)\\[2mm]

\noindent
{
\small 
Dan Abramovich\footnote{Partially supported by NSF grant
DMS-9503276 and by an Alfred P. Sloan research fellowship.} \\
Department of Mathematics, Boston University\\
111 Cummington, Boston, MA 02215, USA\\
{\tt abrmovic@@math.bu.edu}\\[5mm]
Jianhua Wang \\
Department of Mathematics, Massachusetts Institute of Technology\\
77 Massachusetts Avenue,Cambridge, MA 02139, USA\\
{\tt wjh@@math.mit.edu}
}

\addtocounter{section}{-1}
\section{Introduction}

 We work over an algebraically closed field $k$ of characteristic 0.

\subsection{Statement}
In this paper, we use techniques of toric geometry
 to reprove the following theorem: 

\begin{th} Let $X$ be a projective variety of finite type over $k$, and let
$Z\subset X$ be a proper closed subset. Let 
$G\subset \Aut_k(Z\subset X)$ be a finite group. Then there is a
$G$-equivariant 
modification $r:X_1\to X$ such that $X_1$ is nonsingular projective variety,
and $r^{-1}(Z_{\red})$ is a $G$-strict divisor of normal crossings.
\end{th}

This theorem is a weak version of the equivariant case of Hironaka's
well known theorem on resolution of singularities. It was announced by
Hironaka, but a complete proof was not easily accessible for a long time. The
situation was remedied by E. Bierstone and P. Milman \cite{bm1}, who gave a
construction of completely canonical resolution of singularities. Their
construction builds on a thorough understanding of the effect of blowing
up. They carefully build up an invariant pointing to the next blowup. 

The proof we give in this paper takes a completely different approach. It uses
two ingredients: first, we assume that we 
know the existence of resolution of singularities without group actions. The
method of resolution is not important: any of \cite{hi}, \cite{bm}, \cite{aj}
or 
\cite{bp} would do. Second, we use equivariant toroidal resolution of
singularities. Unfortunately, in \cite{te} the authors do not treat the
equivariang case. But proving this turns out to be straightforward given the
methods of \cite{te}.

To this end, section \ref{equitor} of this paper is devoted to proving the
following: 

\begin{th}\label{equitorres}
Let $U\subset X$ be a strict toroidal embedding, and let $G\subset
\Aut(U\subset X)$ be a finite group acting toroidally. Then there is a
$G$-equivariant toroidal ideal sheaf $\I$ such that the normalized blowup of
$X$ along $\I$ is a nonsingular $G$-strict toroidal embedding. 
\end{th}

\subsection{Acknowledgements}
Thanks are due to A. J. de Jong, who was a source of inspiration for this
paper. Thanks also to  S. Katz and T. Pantev for helpful  
discussions relevant to this paper. Special thanks to S. Kleiman, who made this
collaboraion possible. 

\section{Preliminaries}
First recall some definitions. We restrict ourselves to the case of
varieties over $k$. A large portion of the terminology is borrowed from
\cite{aj}. 

A {\bf modification} is a proper birational morphism of irreducible varieties.

Let a finite group $G$ act on a (possibly reducible) variety $Z$. Let $Z=\cup
Z_i$ be the decomposition of $Z$ into irreducible components. We say that
{\bf $Z$ is $G$-strict} if the union of translates $\cup_{g\in G} g(Z_i)$ of
each component $Z_i$ is 
a normal variety. We simply say that $Z$ is
{\bf strict} if it is $G$-strict for the trivial group, namely every $Z_i$ is
normal. 

A divisor $D\subset X$ is called a {\bf divisor of normal crossings} if \'etale
locally at every point it is the zero set of $u_1\cdots u_k$ where
$u_1,\ldots,u_k$ is part of a regular system of parameters. Thus, in a strict
divisor of normal crossings $D$, all components of $D$ are nonsigular.

An open embedding $U\hookrightarrow X$ is called a {\bf toroidal embedding} if
locally in the \'etale topology (or classical topology in case $k=\bfc$, 
or formally) it is
isomorphic to a torus embedding $T \hookrightarrow V$, (see \cite{te}, II\S
1). Let $E_i, i\in I$ be the irreducible components of $X^{_\setminus} U$.  A
finite group action  
$G\subset \Aut(U\hookrightarrow X)$ is said to be {\bf toroidal} if the
stabilizer of every point can be identified on the appropriate neighborhood
with a 
subgroup of the torus $T$. We say that a toroidal action is {\bf $G$-strict} if
$X^{_\setminus} U$ is $G$-strict. In particular the toroidal embedding itself
is 
said to be strict if $X^{_\setminus} U$ is strict. This is the same as the
notion 
of {\bf  toroidal embedding without self-intersections} in \cite{te}. For
any 
subset $J$ of $I$, the components of the sets ${\cap_{i\in J}}E_i -
{\cup_{i\notin J}}E_i$ define a stratification of $X$. Each component
is called a {\bf stratum}.

Recall that in \cite{te}, p. 69-70 one defines the notion of a {\bf conical
polyhedral complex} with {\bf integral structure}. As in \cite{te}, p. 71,
to every strict toroidal embedding $U\subset X$  one canonically
associates a conical polyhedral complex with integral structure. In the sequel,
when we refer to a  conical polhedral complex, it is understood that it is
endowed with an integral structure. 

In \cite{te}, p. 86 (Definition 2) one defines
{\bf a rational finite  partial polyhedral decomposition} $\Delta'$ of a
conical polyhedral complex $\Delta$. We will restrict attention to 
the case where  $|\Delta'| = |\Delta|$, and we will call this simply a {\bf
polyhedral decomposition} or {\bf subdivision}. 

The utility of polyhedral decompositions is given in Theorem 6* of \cite{te}
(page 90), which establishes a correspondence between allowable modifications
of a given strict toroidal embedding (which in our terminology are proper), and
polyhedral decompositions of the conical polyhedral complex. 

In order to guarantee that a modification is projective, one needs a bit more.
Following \cite{te}, p. 91, a function $\ord:\Delta\rar \bfr$ defined on a
conical polydral 
complex with integral structure is called an {\bf order function} if: 

 \parbox{4in}{   (1)   $\ord(\lambda x)=\lambda\cdot{\ord(x)} , \lambda
		\in \bfr^+ $\\
         (2)  $\ord$  is continuous, piecewise-linear\\
         (3) $\ord (N^Y\cap {\sigma^Y}) \subset \bfz     $ for all strata
		$Y$.\\ 
         (4) $\ord$  is convex on each cone $\sigma\subset \Delta$} $(*)$

For an order function on the conical polyhedral complex coresponding to $X$,
we 
can define canonically a coherent sheaf of fractional ideals on $X$, and vice
versa (see
\cite{te}, I\S 2). The order function is positive if and only if the
corresponding sheaf is a  a genuine
ideal sheaf. We have the following important theorem \cite{te}: 

\begin{th} 
Let $F$ be a coherent sheaf of  ideals corresponding to a positive order
function 
$\ord$, and let $B_{F}(X)$ be the normalized blowup of $X$ along $F$. Then
$B_F(X)\rar X$ is
an allowable modification of $X$, described 
by the decompostion of $|\Delta|$ obtained by subdividing  the
cones into the biggest subcones on which $\ord$ is linear. 
\end{th}

A polyhedral decomposition is said to be {\bf projective} if it is obtained in
such a way from an order function.

Given a cone $\sigma$ and a rational ray $\tau\subset \sigma$, it is natural to
defing the decomposition of $\sigma$ centered at $\tau$, whose cones are of the
form $\sigma'+\tau$, where $\sigma'$ runs over faces of $\sigma$ disjoint from
$\tau$. Given a polyhedral complex $\Delta$ and a rational ray $\tau$, we can
take 
the subdivision of all cones containing $\tau$ centered at $\tau$, and again
call the 
resulting decompositionion of $\Delta$, the subdivision centered at $\tau$.

From \cite{te} I\S 2, lemmas 1-3, it follows that the subdivision centered at
$\tau$ is projective. 

A very important decomposition is the {barycentric subdivision}. Let $\sigma$
be a cone with integral structure, $e_1,\ldots,e_k$ integral generators of its
edges. The {\bf barycenter} of $\sigma$ is the ray $b(\sigma) = \bfr_{\geq
0}\sum e_i$. The 
{\bf barycentric subdivision} of a polyhedral complex $\Delta$ of dimension
$m$ is the 
minimal subdivision $B(\Delta)$ in which the barycenters of all cones in
$\Delta$ appear as  cones in $B(\Delta)$. It may be obtained by first taking
the subdivision centered at the barycenters of $m$ dimensional cones, then
the decomposition of the resulting complex centered at the barycenters of the
cones of dimension $m-1$ of the 
{\em original} complex $\Delta$, and so on. From the discussion above (or
\cite{te} III \S 2 lemma 2.2), we have that the barycentric subdivision is 
projective.

One can also obtain the barycentric subdivision inductively the other way: the
barycentric 
subdivision of an $m$-dimensional cone $\delta$  is
formed by first taking the 
barycentric subdivision of all its faces, and for each one of the resulting
cones $\sigma$, including also the cone $\sigma + b(\delta)$. This way it is
clear that 
$B(\Delta)$ is a simplicial subdivision.

\section{Equivariant toroidal modifications}\label{equitor}

\begin{lem}\label{equiact}
   Let $U\subset X$ be a strict toroidal embedding, $G\subset Aut(U\subset X)$
a finite group action. Then
\begin{enumerate}
\item The group $G$ acts linearly on $\Delta(X)$.
\item
 Assume that the action of $G$ is strict toroidal. Let $g\in G$, and let   $
\delta\subset \Delta(X)$ be a cone, such that $ g(\delta)= \delta$.
Then $ g_{|\delta}=id $.
\end{enumerate}
\end{lem}

{\bf Proof.}
\begin{enumerate} \item Clearly, $G$ acts on the stratification of $U\subset
X$. Note that, from Definition 3 of \cite{te}, page 59, 
$\Delta(X)$  is built up from the groups $M^Y$ of Cartier divisors on $Star(Y)$
supported on $Star(Y)^{_\setminus} U$, as $Y$ runs through the strata. As $g\in
G$ 
canonically transforms $M^Y$ to $M^{g^{-1}Y}$ linearly, our claim follows.

\item  Assume $ g: \delta \rightarrow \delta$, and $g_{|\delta}\neq id $, then
there 
 exists an edge $e_1 \in \delta$ , s.t $g(e_1)\neq e_1$. Denote $g(e_1)
=e_2$. Assume $e_1$ corresponds to a divisor $E_1$, and
$e_2$ corresponds to a divisor $E_2$. Since $g(e_1)
=e_2$ we have $g(E_1)=E_2$. As $e_1, e_2$ are both edges of $\delta$,
$E_1\cap E_2  \neq \phi$. So $\cup g(E_1)$ can not be normal since it has two
intersecting components. This is a contradiction to the fact that $G$ acts
strictly on $X$.\qed 
\end{enumerate} 

\begin{lem} \label{equideco}
   Let $G\subset\Aut(U\subset X)$ act toroidally. Let  $\Delta_1$ be a
$G$-equivariant  subdivision of $\Delta$, with corresponding modification
$X_1\rar X$. Then $G$ acts toroidally on $X_1$. Moreover, if $G$ acts strictly
on 
$X$,
it also acts strictly on $X_1$.
\end{lem}

{\bf Proof.}
The fact that $G$ acts on $X_1$ follows from the canonical manner in which
$X_1$ is costructed from the decomposition $\Delta_1$, see
Theorems 6* and 7* of \cite{te}, \S 2.2.
%
%
%

 Now for any point $a \in X_1$ and $g \in Stab_a$, we have $g\circ f(a)=f\circ
g(a)=f(a)$ hence $g \in Stab_{f(a)}$, Thus $Stab_a$ is a subgroup of
$Stab_{f(a)}$,  which is identified with a subgroup of
torus in a neigbourhood of $f(a)$. This proved that $Stab_a$ is identified with
a subgroup of torus.
 
 We are left with showing that if  $G$ acts strictly on $X$, then it acts
strictly on $X_1$.
Assume it is  not the case. There exists two edges $\tau_1$,$\tau_2 $  in
$\Delta_1$, which 
are 
both edges of a cone, $\delta'$, and $g(\tau_1)=\tau_2$. We choose the cone
$\delta'$ of minimal dimension. Clearly, $\tau_1$ 
and $\tau_2 $ cannot  be both edges in $\Delta $, since $G$  acts strictly
on 
$X$. Let us assume $\tau_2$ is not an edge in $\Delta$. So $\tau_2$ must be
in the 
interior of a cone $\delta$ in $\Delta$, which contains $\delta'$. 
Now since ${\delta'}\cap g({\delta'}) \supset \tau_2 \subset$ interior of
$\delta$, 
we conclude: interior of $\delta \cap g(\delta) \neq \phi$, which means
that 
$g(\delta)=\delta$. From the previous lemma, $g_{|\delta}=id$, so 
${g_{|\delta'}}=id$ too, contradiction. \qed

\begin{prp}
\begin{enumerate}
 \item  There is a 1 to 1 correspondence between edges $\tau_i$ in the
barycentric subdivision $B(\Delta)$
and 
positive 
dimensional cones $\delta_i$ in $\Delta $. We denote this by $\tau\mapsto
\delta_\tau$. 
\item Let $\tau_i\neq \tau_j$ be edges of a cone $\hat{\delta}  \in
B(\Delta)$. Then dim $ \delta_{\tau_i}\neq \delta_{\tau_j}. $
\item If $G$ is a finite group acting  toroidally on a strict toroidal
embedding 
$U\subset X$, then  the action of $G$ on 
$  
X_{B(\Delta)}$ is strict. 
\end{enumerate}
\end{prp}

{\bf Remark.} Using this proposition, the argument at the end of \cite{aj} can
be significantly simplified: there is no need to show $G$-strictness of the
toroidal embedding obtained there, since the barycentric subdivision
automatically gives a $G$ strict modification.

{\bf Proof.}
 1. Define a map $b:$ positive dimension cones in $\Delta \rightarrow$ edges in
   $B(\Delta)$ 
by
   $$  b(\delta)=\mbox{ the barycenter of }(\delta)   $$     
      and define $\delta:$ edges in $B(\Delta) \rightarrow $ cones in $\Delta$
      by 
  $$   \delta_\tau= \mbox{ the unique cone whose interior contains }\tau  $$
  then it is easy to see that  $b$ and $\delta$ are invereses of each other.

  2. We proceed by induction on $\dim \Delta$.
The cone $\delta$ spanned by $\tau_i$ and $\tau_j$ must lie in some
cone of 
$\Delta$, say $\delta^*$, which we may take of minimal dimension.
We follow the second construction of the barycentric subdivision described in
the 
 preliminaries. Either $\dim \delta^* \leq m-1$, so
 $\delta$ is in the
barycentric subdivision of 
the $m-1$-skeleton of $\Delta$, in which case the statement follows by the
inductive assumption, or $\dim \delta^*=m$, in which case only one of $\tau_1$
and $\tau_2$ can be its barycenter, and the other is again a barycenter of a
cone in th $m-1$ skeleton.

  3. From lemma \ref{equideco}, since the decomposition $B(\Delta)$ of $\Delta$
is 
equivariant, $G$ acts toroidally on $X_B(\Delta)$. Let $E_1, E_2\subset
X_B(\Delta)^{_\setminus} U$ be divisors
corresponding to edges $ e_1, e_2 $  in $B(\Delta)$. Since $E_1 \cap E_2
\neq 
\phi$, there there is a cone in $B(\Delta)$ containing $e_1, e_2$ as edges.
From part (2),  $\dim \delta_{e_1} \neq  \dim \delta_{e_2}$, so $ g(e_1)$ can
not 
equal 
to $e_2$. This contradicts the fact that the morphism is equivariant
and $g(E_1)=E_2$.

\begin{prp}
   There is a positive $G$-equivariant order function on $B(\Delta)$
such that the associated ideal  $\I$ induces a blowing up
$B_{\I}X_{B(\Delta)}$, which is a nonsingular $G$-strict toroidal embedding, on
which $G$ acts toroidally.
\end{prp}

{\bf Proof.}
 By the previous proposition, we know that  $G$ acts toridally and strictly  on
 $X_{B(\Delta)}$. It follows from Lemma \ref{equiact} that the quotient
   $B(\Delta)/{G}$  is  a conical polyhedral
complex, since no cone has two edges in $B(\Delta)$ which are identified in the
 quotient. We can use the argument 
of \cite{te}, I\S 2, lemmas 1-3, to get an order function $\ord:B(\Delta)/{G}
 \to \bfr$ 
which induces a simplicial  
subdivision with every cell of index 1. Denote by $\pi:B(\Delta) \to
 B(\Delta)/{G}$  the quotient map. Then $ord \circ \pi $ is an order function
 subdividing
$B(\Delta)$ into simplicial cones of  index 1. Let $\I$ be the corresponding
 ideal sheaf. The blow up $X_{B(\Delta)}$
 along $\I$ is  a nonsingular strict toroidal embedding $U\subset
 B_{\I}X_{B(\Delta)}$. By  lemma \ref{equideco}, $G$ acts on
 $B_{\Gamma}X_{B(\Delta)}$ 
strictly and toroidally. \qed

{\bf Proof of Theorem \ref{equitorres}.} Let $G\subset \Aut(U\subset X)$ be as
in the theorem. The morphism $X_{B(\Delta)}\rar X$ is projective, and by the
last two propositions there is a projective, totoidal $G$-equivariant  morphism
$Y \rar X$ where $Y$ is nonsingular and such that $G$ acts strictly and
toroidally on $Y$. \qed

{\bf Remark.} With a little more work we can obtain a {\bf canonical} choice
of a toroidal equivariant resolution of singularities. One observes that
the cones in the barycentric subdivision have canonically ordered coordinates,
which 
agree on intersecting cones: for a cone $\delta$ choose the unit coordinate
vectors $e_i$ to be primitive lattice vectors generating the edges $\tau$,
where $i=\dim \delta_\tau$, the dimension of the cone of which $\tau$ is a
barycenter. Recall that in order to resolve singularities, one successively
takes the 
subdivisions centered at lattice points $w_j$ which are not integrally
generated by the vectors 
$e_i$. These $w_j$ are partially ordered according to the lexicographic
ordering of their canonical coordinates, in such a way that if $w_j\neq w_k$
have the same coordinates (e.g. if $g(w_1) = w_2$), they do not lie in a the
same cone, and therefore we 
can take the centered subdivision simultaneousely.

We conclude this section with a simple proposition which is implicitly used in
\cite{aj}:

\begin{prp} Let $U\subset X$ be a strict toroidal embedding, and let $G\subset
\Aut(U\subset X)$ be a finite group acting strictly and toroidally. Then $
(X/{G},U/{G})$ is a strict toroidal embedding.
\end{prp}

{\bf Proof.} Since the quotient of a toric variety by a finite subgroup of the
torus is toric,  we conclude that $X/{G}$ is still a 
toroidal embedding, by the definition of  toroidal embedding. We 
need to show that it is strict. Let $q:X\rar X/G$ be the quotient map. Let
$Z \subset X^{_\setminus} U$ be a divisor. Then $q(Z) = q(\cup_g g(Z))$. Since
the action is strict, we have $q(\cup_g g(Z)) \simeq Z/Stab(Z)$, which is
normal. 

\section{Proof of the theorem}

Given $Z,X$ with $G$ action , $G$ finite, let $Y=X/G$, $Z/G$ be the
quotient, $B$  the
branch locus. Define $ W = Z/G \cup B$. Let $(Y',W')\to (Y,W)$  be a resolution
of singularities of $Y $
with $W'$ a strict divisor of normal crossings. Let $X'$ be the normalization
of $Y'$ in $K(X)$, 
and $Z'$  the
inverse image of $W'$. Let $U={X'}^{_\setminus} Z'$. Clearly $ U \subset X')$
is 
a strict toroidal embedding, on which
$G$ acts toroidally (moreover, it is $G$-strict). Apply theorem 
\ref{equitorres} and obtain a nonsingular strict toroidal embedding $U\subset
X_1 \rar X'$ as required. \qed

 \end{document}